\documentclass[11pt]{article}
\usepackage{amsfonts,latexsym}
\newcommand{\beq}{\begin{equation}}
\newcommand{\eeq}{\end{equation}}
\newcommand{\beqa}{\begin{eqnarray}}
\newcommand{\eeqa}{\end{eqnarray}}
\newcommand{\beaa}{\begin{eqnarray*}}
\newcommand{\eaa}{\end{eqnarray*}}
\newcommand{\nn}{\hfill\nonumber}
\newcommand \nc {\newcommand}
\nc \proof {{\em{Proof.\/}} }
\nc \qed { $\Box$\hfill}
\newtheorem{theorem}{Theorem}
\newtheorem{lemma}{Lemma}
\newtheorem{proposition}{Proposition}
\newtheorem{corollary}{Corollary}
\newtheorem{definition}{Definition}
\newtheorem{example}{Example}
\nc \bth[1] { \begin{theorem}\label{t#1} }
\nc \ble[1] { \begin{lemma}\label{l#1} }
\nc \bpr[1] { \begin{proposition}\label{p#1} }
\nc \bco[1] { \begin{corollary}\label{c#1} }
\nc \bde[1] { \begin{definition}\label{d#1}\rm }
\nc \bex[1] { \begin{example}\label{e#1}\rm }
\nc \eth { \end{theorem} }
\nc \ele { \end{lemma} }
\nc \epr { \end{proposition} }
\nc \eco { \end{corollary} }
\nc \ede { \end{definition} }
\nc \eex { \end{example} }
\nc \eqref[1] {{\rm{(\ref{#1})}}}
\nc \thref[1]{Theorem \ref{t#1}}
\nc \leref[1]{Lemma \ref{l#1}}
\nc \prref[1]{Proposition \ref{p#1}}
\nc \deref[1]{Definition \ref{d#1}}
\nc \exref[1]{Example \ref{e#1}}
\def \A {{\mathcal A}}

\def \Cset {{\mathbb C}}
\def \Zset {{\mathbb Z}}
\def \Nset {{\mathbb N}}
\def \Vset {{\mathbb V}}
\def \rank { {\mathrm{rank}} }
\def \span { {\mathrm{span}} }
\def \ord { {\mathrm{ord}} }
\def \diag { {\mathrm{diag}} }
\def \spec { {\mathrm{Spec}} }
\def \Im { {\mathrm{Im}} }
\renewcommand \ker { {\mathrm{Ker}} }
\nc \Gr {Gr}
\begin{document}
%%%%%%%%%%%%%%%%%%%%%%%%%%%%%%%%%%%%%%%%%%%%%%%%%%%%%%%%%%%%%%%%%%%%%%%
%%%%%%%%%%%%%%%%%%%%%%    Title    %%%%%%%%%%%%%%%%%%%%%%%%%%%%%%%%%%%%%%%%%%%
\title{{\LARGE\bf{B\"acklund--Darboux transformations in Sato's Grassmannian}}}
\author{
B.~Bakalov
\thanks{E-mail: bbakalov@fmi.uni-sofia.bg}
\quad
E.~Horozov
\thanks{E-mail: horozov@fmi.uni-sofia.bg}
\quad
M.~Yakimov
\thanks{E-mail: myakimov@fmi.uni-sofia.bg}
\\ \hfill\\ \normalsize \textit{
Department of Mathematics and Informatics, }\\
\normalsize \textit{Sofia University , 5 J. Bourchier Blvd.,
Sofia 1126, Bulgaria }     }
\date{}
\maketitle
\begin{abstract}
%\hspace*{-22pt}
We define B\"acklund--Darboux transformations in Sato's Grassmannian.
They can be regarded as Darboux transformations on maximal
algebras of commuting ordinary differential operators. We describe the
action of these transformations on related objects:
wave functions, tau-functions and spectral algebras.
\end{abstract}
%%%%%%%%%%%%%%%%%%%
\vspace{-11cm}
\begin{flushright}
{\tt{ q-alg/9602010 }}
\end{flushright}
\vspace{9cm}
%%%%%%%%%%%%%%%%%%%
%%%%%%%%%%%%%%%  Sect. 0  %%%%%%%%%%%%%%%%%%%%%%%%%%%%%%%%%%%%%%%%%%%%%%%%%%
\setcounter{section}{-1}
\section{}
%{Introduction}
%%%%%%%%%%%%%%%%%%%%%%%%%%%%%%%%%%%%%%%%%%%%%%%%%%%%%%%%%%%%%%%%%%%%%%%%%%%%

Classically, a Darboux transformation \cite{Da} of a differential operator $L$,
presented as a product $L = Q P$, is defined by exchanging the places of the
factors, i.e.\ $ \overline L = P Q $.
Obviously all versions of Darboux transformations have the property that if
$\Phi(x)$ is an eigenfunction of $L$, i.e.\ $L\Phi=\lambda\Phi$ then $P\Phi$ is
an eigenfunction of $\overline L$, i.e.\ $\overline LP\Phi=\lambda P\Phi$.
Motivated by this characteristic property we give a version of Darboux
transformation directly on wave functions. The plane $W$ of Sato's
Grassmannian  is said to be a {\it{B\"acklund--Darboux transformation}} of the
plane $V$ iff the corresponding wave functions are connected by:
\beqa
&&\Psi_W(x,z)=\frac{1}{g(z)} P(x,\partial_x) \Psi_V(x,z),
\nn\hfill\\
&&\Psi_V(x,z)=\frac{1}{f(z)} Q(x,\partial_x) \Psi_W(x,z)
\nn\hfill
\eeqa
for some polynomials $f$, $g$ and differential operators $P$, $Q$, or
equivalently --
$$fV\subset W\subset \frac1{g} V.$$

To any plane $W\in\Gr$ one can associate a maximal algebra $\A_W$  of
commuting ordinary differential operators \cite{Kr, S, SW} (called a
{\em{spectral
algebra\/}}). Recall that a {\em{rank\/}} of $\A_W$ is the g.c.d. of the orders
of the operators from $\A_W$. We prove that B\"acklund--Darboux
transformations preserve the rank of the spectral algebra. Moreover if $W$ is
a  B\"acklund--Darboux transformation of $V$ such that $\A_V=\Cset[L_V]$
for some operator $L_V$ then every operator from $\A_W$ is a Darboux
transformation (in the sense of eq.\ \eqref{3.7} below) of an operator from
$\A_V$.

In our terminology the set of rational solutions of the KP hierarchy \cite{Kr,
SW} coincides with the set of B\"acklund--Darboux transformations of the
simplest plane
$H_{+} = \span\{ z^n | n=0,1,\ldots\}$. The corresponding tau-function $\tau_W$
is given by the so-called ``{\it{superposition low for wobbly solitons\/}}''
(see e.g.\ \cite{SW, W}). We generalize this formula for a B\"acklund--Darboux
transformation $W$ of an arbitrary plane $V$ provided that
the wave function $\Psi_V(x,z)$ is well defined at all zeros of the
polynomial $f(z) g(z)$ (see \thref{3.12} below). In a particular but important
case which we use it is proven in \cite{AvM2}. The case when
$\Psi_V(x,z)$ is not well defined for some zero $z = \lambda$ of $f(z) g(z)$
is even more interesting (cf.\ \cite{BHY2, BHY3}). We obtain a formula for
$\tau_W$ valid in the general situation (see \thref{3.8} below).

A geometric interpretation of $\ker P$ can be given using the so-called
conditions $C$ (introduced in \cite{W} for the rational solutions of
KP). 
When the spectral curve $\spec\A_V=\Cset$ (i.e.\ $\A_V=\Cset[L_V]$)
the spectral curve $\spec \A_W$ of $W$ can be obtained from that of $V$ by
introducing singularities at points where the conditions $C$ are supported
-- see \cite{W} and also \cite{BHY2} (from \cite{BC} it is known that $\spec
\A_W$ is an algebraic curve).

This paper may be considered as a part of our project \cite{BHY1}--\cite{BHY3}
on the bispectral problem (see \cite{DG}). Although here we do not touch the
latter, the present paper arose in the process of working on \cite{BHY2,BHY3}.
We noticed that many facts, needed in \cite{BHY2,BHY3} can be naturally
obtained in a more general situation. Apart from the applications to the
bispectral problem we hope that some of the results can be useful elsewhere.

{\flushleft{\bf{Aknowledgement}}}

\medskip\noindent
This work was partially supported by Grant MM--523/95 of Bulgarian
Ministry of Education, Science and Technologies.

%%%%%%%%%%%%%%%  Sect. 1   %%%%%%%%%%%%%%%%%%%%%%%%%%%%%%%%%%%%%%%%%%%%%%%%%%
\section{}
%{Preliminaries on Sato's Grassmannian }
%%%%%%%%%%%%%%%%%%%%%%%%%%%%%%%%%%%%%%%%%%%%%%%%%%%%%%%%%%%%%%%%%%%%%%%%%%%%

The aim of this section is to recollect some facts and notation from Sato's
theory of KP-hierarchy \cite{S, DJKM, SW} needed in the paper.
We  use the approach of  V.\ Kac and D.\ Peterson  based on infinite wedge
products (see e.g.\ \cite{KRa}) and the recent survey paper by P. van~Moerbeke
\cite{vM}.

Consider the space of formal Laurent series in $z^{-1}$
$$
\Vset=\Bigl\{ \sum_{k\in\Zset} a_kz^k \Big|\; a_k=0\ {\rm for}\
k\gg 0\Bigr\}.
$$
We define the fermionic  Fock space $F$ consisting of formal infinite sums of
semi-infinite monomials
$$
z^{i_0}\wedge z^{i_1}\wedge\ldots
$$
such that $i_0<i_1<\ldots$ and $i_k=k$ for $k\gg0$.
Let $gl_\infty$ be the Lie algebra of all $\Zset\times \Zset$ matrices having
only a finite number of non-zero entries. One can define a representation $r$
of $gl_\infty$ in the fermionic Fock space $F$ as follows
\beq
r(A)(z^{i_0}\wedge z^{i_1}\wedge\ldots) =
       Az^{i_0}\wedge z^{i_1}\wedge\ldots +
       z^{i_0}\wedge Az^{i_1}\wedge\ldots + \cdots. \label{1.1'9}
\eeq
 The above defined representation $r$ obviously cannot be continued on
the Lie algebra $\widetilde{gl}_{\infty}$ of all $\Zset\times \Zset$ matrices
with finite number of non-zero diagonals.
If we regularize it by
\beq
\hat r(D)(z^{i_0}\wedge z^{i_1}\wedge\ldots) =
              \bigr[(d_{i_0}+d_{i_{1}}+\cdots)-
                    (d_0+d_1+\cdots) \bigl]
              (z^{i_0}\wedge z^{i_1}\wedge\ldots)
\label{1.2}
\eeq
for $D = \diag(\ldots,d_{-1},d_{0},d_{1},\ldots)$ and by \eqref{1.1'9} for an
off-diagonal matrix this will give a representation of a
central extension $\widehat{gl}_\infty =\widetilde{gl}_\infty \oplus \Cset c$
of $\widetilde{gl}_\infty$. Here the central charge $c$ acts as multiplication
by 1. Introduce also the shift matrices $\Lambda_n(n \in \Zset)$ representing
the multiplication by $z^n$ in the basis $\{ z^i\}_{i \in \Zset}$ of $\Vset$.
Then $\hat r(\Lambda_k)$ generate a representation of the Heisenberg algebra:
$$ \Bigl[\hat r(\Lambda_n),\hat r(\Lambda_m)\Bigr] = n \delta_{n,-m}. $$
There exists a unique isomorphism (see \cite{KRa} for details):
\begin{eqnarray}
&&\sigma\colon F\to B=\Cset\left[[t_1, t_2, t_3,\dots]\right]
\hfill\label{2,5}\\
&&\sigma(\hat r(\Lambda_n))= \frac{\partial}{\partial t_n},\quad
\sigma(\hat r(\Lambda_{-n}))=nt_n,\qquad n>0,\quad   \hfill \label{1.4'6}
\end{eqnarray}
known as the boson-fermion correspondence ($B$ is called a bosonic Fock space).

{\em Sato's Grassmannian\/} $\Gr$ \cite{S, DJKM, SW} consists of all
subspaces $W\subset\Vset$ which have an admissible basis
$$
w_k=z^k+\sum_{i<k}w_{ik}z^i,\quad k=0,1,2,\ldots
$$
To a plane $W\in \Gr$ we associate a state
$|W\rangle\in F$ as follows
$$
|W\rangle = w_0\wedge w_{1}\wedge w_{2}\wedge\ldots
$$
A change of the admissible basis results in a multiplication of $|W\rangle$ by
a non-zero constant. Thus we define an embedding of $\Gr$ into the
projectivization of $F$ which is called a Pl\"ucker
embedding.
One of the main objects of Sato's theory is the {\em  tau-function\/} of $W$
defined as the image of $|W\rangle$ under the boson-fermion correspondence
\eqref{2,5}:
\beq
\tau_W(t) =\sigma(|W\rangle) =\langle 0|\; e^{H(t)}\;|W\rangle,
\label{1.6}
\eeq
where $H(t)=-\sum_{k=0}^\infty t_k {\hat r}(\Lambda_k)$.
Another important function connected to $W$ is the {\em  Baker\/} or {\em wave
function\/}
\beq
\Psi_W(t,z)= e^{\sum_{k=1}^\infty t_kz^k}
\frac{\tau\left(t-[z^{-1}]\right)}{\tau(t)},
\label{1.7}
\eeq
where $[z^{-1}]$ is the vector $\left(z^{-1}, z^{-2}/2,\ldots\right)$.
We often use the notation
$\Psi_W(x,z)$ $=\Psi_W(t,z)|_{t_1=x, t_2=t_3=\cdots=0}$.

The Baker function $\Psi_W(x,z)$ contains the whole information about $W$ and
hence about $\tau_W$, as the vectors
$w_{k}=\partial^k_x\Psi_W(x,z)|_{x=0}$
form an admissible basis of $W$. We can expand $\Psi_W(t,z)$ in a formal
series as
\beq
\Psi_W(t,z)=e^{\sum_{k=1}^\infty t_kz^k} \left(1+\sum_{k>0}a_k(t)z^{-k}\right).
\label{1.13}
\eeq
Introduce also the  pseudo-differential operators
$K_W(t,\partial_x) =1+\sum_{j>0}a_j(t)\partial_x^{-j}$ such that
$\Psi_W(t,z) = K_W(t,\partial_x) e^{\sum_{k=1}^\infty t_k z^k}$
(the {\em{wave
operator\/}}) and $P=K_W \partial_x K_W^{-1}$.
Then $P$ satisfies the following infinite system of non-linear differential
equations
\beq
\frac{\partial}{\partial t_k} P=\Bigl[P^k_+, P\Bigr],
\label{1.16}
\eeq
called the {\em KP hierarchy} and $\Psi_W(x,z)$ is an eigenfunction of
$P(x,\partial_x)$:
\beq
P\Psi_W(x,z) = z\Psi_W(x,z).
\label{1.17}
\eeq
A very important  object connected to the plane $W$ is the
algebra $A_W$ of polynomials $f(z)$ that leave $W$ invariant:
\beq
A_W=\{ f(z) | f(z)W\subset W\}.
\label{1.20}
\eeq
For each $f(z)\in A_W$ one can show that there exists a unique
differential operator $L_f(x,\partial_x)$, the order of $L_f$ being equal to
the degree of $f$, such that
\beq
L_f\Psi_W(x,z)=f(z)\Psi_W(x,z).
\label{1.20'5}
\eeq
Explicitly we have
\beq
L_f=K_Wf(\partial_x)K_W^{-1}.
\label{1.21}
\eeq
We denote the commutative algebra of these operators by $\A_W$, i.e.\
\beq
\A_{W} = \{ L_{f} | L_{f} \Psi_{W} = f \Psi_{W}, \; f\in A_{W} \}.
\label{1.22}
\eeq
Obviously $A_W$ and $\A_W$ are isomorphic.
We call $\A_W$ {\em spectral algebra\/} corresponding to the plane $W$.
Following I.~Krichever
\cite{Kr} we introduce the {\em  rank\/} of $\A_W$ to be the dimension of
the space of joint eigenfunctions of the operators from $\A_W$. It coincides
with the greatest common divisor of the orders of the operators $L_f$.

We also use the notation $\Gr^{(N)} := \{W \in \Gr | z^N \in A_W \}$. It
coincides with the subgrassmannian of solutions of the so-called $N$-th
reduction of the KP hierarchy.
%%%%%%%%%%%%%%%%%%  Sect. 2             %%%%%%%%%%%%%%%%%%%%%%%%%%
%%%%%%%%%%%%%%%%%%% % Definitions                        %%%%%%%%%
%%%%%%%%%%%%%%%%%%%%%%%%%%%%%%%%%%%%%%%%%%%%%%%%%%%%%%%%%%%%%%%%%%
\section{}

In this section we introduce our basic definition of {\em B\"acklund--Darboux
transformation} in the Sato's Grassmannian.

The classical Darboux transformation is defined on ordinary
differential operators in the variable $x$, presented in a factorized form
$L=Q P$; it exchanges the places of the factors, i.e.\ the image of $L$ is
the operator $\overline L=P Q$.
The next classical lemma answers the question when the factorization
$L=QP$ is possible (see e.g.\ \cite{I}).
\ble{3.1}
$L$ can be factorized as
\beq
L=QP\ \ {\it iff}\ \ \ker P\subset \ker L.
\label{3.4}
\eeq
In this case
\beq
\ker Q = P(\ker L).
\label{3.5}
\eeq
\ele
A slightly more general construction is the following one.
For operators $L$ and $P$ such that the kernel of $P$ is invariant under $L$,
i.e.\
\beq
L(\ker P)\subset \ker P
\label{3.6}
\eeq
we consider the transformation
\beq
L\mapsto \overline L=PLP^{-1}.
\label{3.7}
\eeq
The fact that $\overline L$ is a differential operator follows from
\leref{3.1}. Indeed, $L(\ker P)\subset \ker P$ is equivalent to $\ker P\subset
\ker(PL)$.
If $h$ is the characteristic polynomial of the linear operator $L|_{\ker
P}$ then $h(L|_{\ker P})=0$, i.e.\
$$
\ker P\subset \ker h(L).
$$
This shows that
$$
h(L)=QP,\ \ h(\overline L)=PQ
$$
for some operator $Q$.

Now we come to our basic definition.
\bde{3.2}
We say that a plane $W$ (or the corresponding wave function $\Psi_W(x,z)$) is a
{\em B\"acklund--Darboux transformation\/} of the plane $V$ (respectively wave
function $\Psi_V(x,z)$) iff there exist (monic) polynomials $f(z)$, $g(z)$ and
differential operators $P(x,\partial_x)$, $Q(x,\partial_x)$ such that
\beq
\Psi_W(x,z)=\frac{1}{g(z)} P(x,\partial_x) \Psi_V(x,z),
\label{3.8}
\eeq
\beq
\Psi_V(x,z)=\frac{1}{f(z)} Q(x,\partial_x) \Psi_W(x,z).
\label{3.9}
\eeq
\ede

Here necessarily $\ord P= \deg g$ and $\ord Q= \deg f$. 
The polynomial $g(z)$ can be chosen arbitrary but of the same degree
because the wave function $\Psi_W(x,z)$ is determined up to a
multiplication by a formal series of the form
$
1 + \sum_{k=1}^\infty a_k z^{-k}.
$

Note that a
composition of two B\"acklund--Darboux transformations is again a
B\"acklund--Darboux transformation. For example the bispectral potentials of
\cite{DG} can be obtained by one B\"acklund--Darboux transformation in
contrast to the finite number of ``rational'' Darboux transformations of
Duistermaat and Gr\"unbaum.

Simple consequences of \deref{3.2} are the identities
\beq
PQ\Psi_W(x,z)=f(z)g(z)\Psi_W(x,z),
\label{3.10}
\eeq
\beq
QP\Psi_V(x,z)=f(z)g(z)\Psi_V(x,z),
\label{3.11}
\eeq
i.e.\ the operator $\overline L=PQ \in \A_W$ is a Darboux transformation of
$L=QP \in \A_V.$ Obviously \eqref{3.8} implies the inclusion
\beq
gW\subset V.
\label{3.12}
\eeq
Conversely, if \eqref{3.12} holds there exists $P$ satisfying \eqref{3.8}.
Therefore a definition equivalent to \deref{3.2} is the following one.
\bde {2}
A plane $W$ is a B\"acklund--Darboux transformation of a plane $V$ iff
\beq
fV\subset W\subset \frac1gV
\label{3.13}
\eeq
for some polynomials $f(z)$, $g(z)$.
\ede
%%%%%%%%%%%%%%%%%%%%%%%%% Sect.3 %%%%%%%%%%%%%%%%%%%%%%%%%%%%%%%%%%%%%%%%
% B\"acklund--Darboux transformations on the spectral algebras
%%%%%%%%%%%%%%%%%%%%%%%%%%%%%%%%%%%%%%%%%%%%%%%%%%%%%%%%%%%%%%%%%%%%%%%%%
\section{}

In this section we study the behavior of the spectral algebra under
B\"acklund--Darboux transformations. The following simple lemma will be useful.
\ble{3.3}
 In the notation of\/ {\em \eqref{3.13}} and\/    {\em \eqref{1.20}}
\beq
fgA_V\subset A_W\subset \frac1{fg} A_V. \label{3.13'6}
\eeq
\ele
\noindent
The {\it proof} is obvious from \eqref{3.13}.
\qed

\smallskip\noindent
\bpr{3.4}
The B\"ackund--Darboux transformations preserve the rank
of the spectral algebras, i.e.\ if $W$ is a B\"ackund--Darboux transformation
of $V$ then $\rank \A_W=\rank \A_V.$
\epr
\noindent
\proof
Let $\rank A_V=N$. Then \leref{3.3} implies that $f g \in A_V$ and
therefore $\deg fg = Nj$, $j\in \Nset$. The right inclusion in
\eqref{3.13'6}
gives $ N | \rank A_W$. Because $\rank A_V=N$, $A_V$ contains polynomials of
degrees $ka+lb$, for $k,l \in \Zset_{\geq 0}$ with $(a,b)=N$. The left
inclusion in \eqref{3.13'6} implies that $A_W$ contains polynomials of degrees
$ka+lb+Nj$, for $k,l \in \Zset_{\geq 0}$, i.e.\ $\rank A_W | N$ and
therefore $\rank A_W=N$.
\qed

\smallskip\noindent
The most important case of algebras $\A_V$ of $\rank N$ is
\beq
A_V=\Cset[z^N],\quad \A_V=\Cset[L_V]
\label{3.14}
\eeq
for some natural number $N$ and a differential operator $L_V$ of order
$N$ (see \cite{BHY2}). This corresponds to the case when the spectral curve
$\spec \A_V$ of $V$ is $\Cset$. We shall describe $\A_W$ for a
B\"acklund--Darboux transformation
$W$ of $V$ for which \eqref{3.14} holds. First observe that due to \eqref{3.11}
we have
\beqa
&&f(z)g(z)=h(z^N),
\label{3.15}   \hfill\\
&&QP=h(L_V)
\label{3.16}   \hfill
\eeqa
for some polynomial $h$.
\bpr{3.4a}
If $\A_V=\Cset[L_V],$  $\ord L_V=N$ then
\beqa
&&A_W=\left\{ u\in \Cset[z^N] \mid u(L_V)\ker P\subset \ker P\right\},
\label{3.18}\hfill\\
&&\A_W=\left\{ Pu(L_V)P^{-1} \mid u\in A_W\right\}.
\label{3.19}
\eeqa
\epr
\noindent
\proof
Since $A_W \subset \Cset[z]$ the right inclusion of \eqref{3.13'6} with
$A_V = \Cset[z^N]$ and $f(z)g(z)=h(z^N)$ implies $A_W \subset \Cset[z^N]$.
Let $u(z^N)\in A_W$
and $L$ be the corresponding operator from $\A_W$ (see \eqref{1.21}), such
that $$
L\Psi_W(x,z)=u(z^N)\Psi_W(x,z).
$$
Using \eqref{3.8} we compute
$$
L\Psi_W(x,z)=L\frac1{g(z)} P\Psi_V(x,z)=\frac 1{g(z)} LP\Psi_V(x,z)
$$
and
\beqa
&&u(z^N)\Psi_W(x,z) = u(z^N) \frac1{g(z)} P\Psi_V(x,z)
\nn\\
&&= \frac1{g(z)} Pu(z^N)\Psi_V(x,z) = \frac1{g(z)} Pu(L_V)\Psi_V(x,z).
\nn
\eeqa
Therefore
$$
LP=Pu(L_V).
$$
The operator $L=Pu(L_V)P^{-1}$ is differential iff $u(L_V)\ker P\subset \ker P$
and obviously it belongs to $\A_W$.
\qed

\smallskip\noindent
Thus the determination of $A_W$ is reduced to the following finite-dimensional
problem:

\smallskip
{\em For the linear operator $L_V|_{\ker h(L_V)}$ find all polynomials
$u(z^N)$ such that the subspace\/ $\ker P$ of\/ $\ker h(L_V)$ is invariant
under the operator $u(L_V)$.}
%%%%%%%%%%%%%%%%%%%%%%%%%%%%%%%%%%%%%%%%%%%%%%%%%%%%%%%%%%%%%%%%%%%%%%%%%%%%%%
% B\"acklund--Darboux transformations on tau-functions
%%%%%%%%%%%%%%%%%%%%%%%%%%%%%%%%%%%%%%%%%%%%%%%%%%%%%%%%%%%%%%%%%%%%%%%%%%%%%%
\section{}

In this section we shall describe the tau-function of the B\"acklund--Darboux
transformation $W$ of $V$ in terms of the tau-function of $V$ and $\ker P$.

Let $V\in \Gr^{(N)}$, i.e.\ $z^NV\subset V$ and
$L_V\Psi_V(x,z)=z^N\Psi_V(x,z)$ for some operator $L_V$ of order $N$.
(We do not suppose that $A_V=\Cset[z^N]$ but only that $\Cset[z^N]\subset
A_V$.) Let $W$ be a B\"acklund--Darboux transformation of $V$ such that
(\ref{3.8}, \ref{3.9}, \ref{3.15}, \ref{3.16}) hold. Let us fix a basis
$\{\Phi_i(x)\}_{0\le i\le dN-1}$ of $\ker h(L_V)$ (where $d=\deg h$). The
kernel of $P$ is a subspace of $\ker h(L_V)$. We fix a basis of $\ker P$
\beq
f_k(x)=\sum_{i=0}^{dN-1} a_{ki}\Phi_i(x),\quad 0\le k\le n-1.
\label{3.24}
\eeq
We can suppose that $P$ and $g$ are monic. Eq.\ \eqref{3.8} implies (see e.g.\
\cite{I})
\beq
\Psi_W(x,z)=\frac{Wr\bigl( f_0(x),\ldots, f_{n-1}(x),\Psi_V(x,z)\bigr)}
{g(z)Wr\bigl( f_0(x),\ldots, f_{n-1}(x)\bigr)},
\label{3.25}
\eeq
where $Wr$ denotes the Wronski determinant.
When we express $f_k(x)$ by \eqref{3.24} we obtain
\beq
\Psi_W(x,z)=
\frac{\sum\det A^I Wr\bigl(\Phi_I(x)\bigr) \Psi_I(x,z)}
{\sum\det A^I Wr\bigl(\Phi_I(x)\bigr)}.
\hfill \label{3.26}
\eeq
The sum is  taken over all $n$-element subsets
$$
I=\{i_0<i_1<\ldots <i_{n-1}\} \subset\{0,1,\ldots, dN-1\}
$$
and here and further we use the following notation:
$$
A=(a_{ki})_{0\le k\le n-1,\; 0\le i \le dN-1}
$$
is the matrix from \eqref{3.24} and
$$
A^I=(a_{k,i_l})_{0\le k\le n-1,\; 0\le l\le n-1}
$$
is the corresponding minor of $A$,
$$
\Phi_I(x)=\left\{\Phi_{i_0}(x),\ldots,\Phi_{i_{n-1}}(x)\right\}
$$
is the corresponding subset of the basis $\{\Phi_i(x)\}$ of $\ker h(L_V)$ and
\beq
\Psi_I(x,z)=\frac{Wr\bigl(\Phi_I(x),\Psi_V(x,z)\bigr)}{g(z)Wr\bigl( \Phi_I(x)
\bigr)}
\label{3.27}
\eeq
is the B\"acklund--Darboux transformation of $V$ with a basis of  $\ker P$
$f_k(x)=\Phi_{i_k}(x)$.

So if we know how $V$ transforms when the basis $\{f_k(x)\}$ of $\ker P$ is a
subset of the basis $\{\Phi_i(x)\}$ of $\ker h(L_V)$, the formula \eqref{3.26}
gives us $\Psi_W(x, z)$ for an arbitrary B\"acklund--Darboux transformation
$W$.

We shall obtain a similar formula for the tau-functions as well.

Let $\tau_V(0)$ and
$\tau_W(0)$ be nonzero and let us normalize them to be equal to 1
(recall that the tau-function is defined up to a multiplication by a constant).
We set
$$
\Delta_I= Wr\bigl(\Phi_I(x)\bigr)|_{x=0}.
$$
Denote by $\tau_I$ the tau-function corresponding to the wave function
\eqref{3.27}, also normalized by $\tau_I(0)=1$.
Then we have the following theorem.
\bth{3.8}
 In the above notation
\beq
\tau_W(t)=\frac{\sum\det A^I\Delta_I\tau_I(t)}{\sum\det A^I\Delta_I}.
\label{3.28}
\eeq
\eth
\noindent
For the {\em proof\/} we have to introduce some more terminology. These are the
so-called {\em conditions\/} $C$ (cf.\ \cite{W}).  They are conditions (or
equations)
that should be imposed on a vector $v\in V$ in order to belong to $gW$ (recall
\eqref{3.12}).

Let us fix an admissible basis $\{ v_k\}_{k\ge0}$ of $V$ and set
$$
V_{(n)}=\bigoplus_{k=0}^{n-1}\Cset v_k
$$
(this is independent of the choice of the basis).
We define a linear map
$$
C\colon V\to V_{(n)}
$$
by defining it on the wave function of $V$
\beq
C\Psi_V(x,z)=\sum_{k=0}^{n-1} f_k(x) v_k,
\label{3.29}
\eeq
where $\{f_k(x)\}$ is the basis of $\ker P$.
The point is that $C$ acts on the variable $z$.
If we choose $v_k=\partial_x^k\Psi_V(x,z)|_{x=0}$ then
$$
Cv_p=C\partial_x^p\Psi_V(x,z)|_{x=0}=
\partial^p_x C\Psi_V(x,z)|_{x=0}=
\sum_{k=0}^{n-1} f_k^{(p)}(0) v_k.
$$
Let $V_C$ be the kernel of $C$, i.e.\
$$
V_C=\{v\in V\mid Cv=0\}.
$$
Then the description of $gW$ is straightforward (cf.\ \cite{W}).
\ble{3.9}
$W={\displaystyle\frac{1}{g}}V_C$.
\ele
\noindent
{\it Proof}. First we show that $gW\subset V_C$. Indeed,
$$
C\bigl(g(z)\Psi_W(x,z)\bigr)= C\bigr(P(x,\partial_x)\Psi_V(x,z)\bigr)=
P(x,\partial_x) C\Psi_V(x,z)
$$
$$
{}=P(x,\partial_x)\sum_{k=0}^{n-1} f_k(x) v_k=0,\quad
{\rm because}\ f_k\in\ker P.
$$
On the other hand the vectors $g(z)\partial_x^j\Psi_W(x,z)|_{x=0}$ can be
expressed in the form
$$
v_{j+n}+\sum_{k<j+n} d_{jk}v_k,\quad j\ge 0,
$$
i.e.\ the plane $gW$ maps one to one on the plane $\bigoplus_{j\ge n}\Cset
v_j$. But the same is true for the plane $V_C$ as $\Im C=V_{(n)}$ (because
$\det(C|_{V_{(n)}})=Wr\bigl( f_k(x)\bigr)|_{x=0}\ne 0$).
\qed
\bco{3.10}
$W$ has an admissible basis
\beq
w_j=\frac1{g(z)}\left(1-C^{-1}_{(n)} C\right) v_{j+n},\quad j\ge 0,
\label{3.30}
\eeq
where $C_{(n)}=C|_{V_{(n)}}$.
\eco
\noindent
{\it Proof}.
$C\left(1-C^{-1}_{(n)}C\right) =C-C_{(n)}C^{-1}_{(n)}C=0$ since $\Im
C=V_{(n)}$.
\qed

\smallskip\noindent
$V_C$ can also be interpreted as the intersection of the kernels of certain
linear functionals on $V$
(${\mathrm{pr}}_k\circ C\colon V\to\Cset v_k\equiv\Cset$)
which form an $n$-dimensional linear
space. We denote it by abuse of notation again by $C$.
\ble{3.11}
 Any condition $c\in C$ gives rise to a function
\beq
f(x)=\langle c,\Psi_V(x,z)\rangle
\label{3.31}
\eeq
{}from\/ $\ker P$, and vice versa.
\ele
\noindent
The {\it proof\/} follows immediately from the definition \eqref{3.29}.
\qed

\smallskip\noindent
We define linear functionals $\chi_i$ and $c_k$ on $V$ by
\beqa
&&\langle \chi_i,\Psi_V(x,z)\rangle =\Phi_i(x),\quad 0\le i\le dN-1,
\label{3.32} \hfill\\
&&\langle c_k,\Psi_V(x,z)\rangle =f_k(x),\quad 0\le k\le n-1,
\label{3.33}
\eeqa
i.e.\ $c_k=\sum_{i=0}^{dN-1} a_{ki}\chi_i$.

\smallskip\noindent
We can now give the {\it {proof of \thref{3.8}}}.
The basis \eqref{3.30} of $W$ can be written as
$$
w_j=\frac1{g(z)} \Bigl( v_{j+n} -\sum_{0\le k\le n-1,\, 0\le i\le dN-1}
\left(C^{-1}_{(n)} A\Bigr)_{ki}
\langle \chi_i,v_{j+n}\rangle v_k\right), \quad j\ge 0.
$$
We use the formula \eqref{1.6}
$$
\tau_W(t)=\sigma(w_0\wedge w_1\wedge w_2 \wedge\ldots)
$$
and expand all $w_{j}$:
$$
\tau_W(t)=\sum_{r=-1}^{n-1}
\sum_{\textstyle{0\le k_s\le n-1,\; 0\le i_s\le dN-1 \atop {\rm for}\ 0\le s\le
r}} \left(C^{-1}_{(n)}A\right)_{k_0i_0}
\cdots
\left(C^{-1}_{(n)}A\right)_{k_ri_r}\cdot(-1)^{r+1}
$$
$$
\times  \sum_{n \le j_0 < \ldots < j_r}
\sigma\left( \frac1g v_n\wedge \frac1g v_{n+1}\wedge \ldots
\wedge \frac1g\langle \chi_{i_0},v_{j_0}\rangle v_{k_0}\wedge\ldots
\wedge \frac1g \langle \chi_{i_r}, v_{j_r}\rangle v_{k_r} \wedge\ldots\right)
$$
(the term
$\frac1g\langle \chi_{i_s},v_{j_s}\rangle v_{k_s}$
is on the ($j_s - n + 1$)-st place in the wedge product).
Let $(k_0,\ldots , k_{n-1})$ be a permutation of $\{0,1,\ldots,n-1\}$.
Noting that
$$
\sum_{0\le i\le dN-1} \left( C^{-1}_{(n)} A\right)_{k i}
\langle \chi_i,v_j\rangle=\delta_{kj} \quad{\rm for}\
0\le k,j\le n-1
$$
we can insert
$$
\sum_{\textstyle{0\le k_s\le n-1,\; 0\le i_s\le dN-1 \atop {\rm for}\ r+1\le
s\le n-1}} \left(C^{-1}_{(n)}A\right)_{k_{r+1}i_{r+1}}
\cdots
\left(C^{-1}_{(n)}A\right)_{k_{n-1}i_{n-1}}
$$
$$
\times  \sum_{0 \le j_{r+1} < \ldots < j_{n-1} \le n-1}
\langle \chi_{i_{r+1}},v_{j_{r+1}}\rangle  \cdots
\langle \chi_{i_{n-1}}, v_{j_{n-1}}\rangle 
$$
in the above expression for $\tau_W(t)$.
Then
$$
\tau_W(t)=\sum_{(k_0,\ldots,k_{n-1})}
\sum_{0\le i_0,\ldots, i_{n-1}\le dN-1}
\left(C^{-1}_{(n)}A\right)_{k_0i_0}
\ldots
\left(C^{-1}_{(n)}A\right)_{k_{n-1}i_{n-1}}
$$
$$
\times \, \sigma\left(
R\left(\frac1g\right) r(\chi_{i_0})\cdots r(\chi_{i_{n-1}})
(v_{k_0}\wedge\ldots \wedge v_{k_{n-1}}\wedge v_n \wedge v_{n+1}\wedge \ldots
)\right),
$$
where the operator $R\bigl(\frac1g\bigr)$ acts as a group element
$$
R\left(\frac1g\right) (u_0\wedge u_1\wedge\ldots)=
\frac1g u_0\wedge \frac1g u_1\wedge \ldots
$$
and $r(\chi_i)$ is a contracting operator:
$$
r(\chi_i)(u_0\wedge u_1\wedge\ldots)=
\sum_{j\ge0} (-1)^j \langle \chi_i, u_j\rangle
u_0\wedge u_1\wedge\ldots \wedge \widehat u_j\wedge\ldots
$$
(the hat on $u_j$ means as usually that it is omitted).

By the antisymmetry we obtain
\beq
\tau_W(t)=\sum\det\left( C_{(n)}^{-1}A\right)^I R\left(\frac1g\right)
r(\chi_I)\tau_V(t),
\label{3.34}
\eeq
where the sum is over the subsets
$$
I=\{ i_0<\ldots <i_{n-1}\} \subset \{0,1,\ldots,dN-1\}
$$
and
$$
r(\chi_I)=r(\chi_{i_0})\cdots r(\chi_{i_{n-1}}).
$$
For the special B\"acklund--Darboux transformation $\tau_I(t)$ with
$f_k(x)=\Phi_{i_k}(x)$ we
have
$$
c_k=\chi_{i_k},\quad A=(\delta_{ii_k})_{\textstyle{0\le k \le n-1\atop 0\le i
\le dN-1}}, \quad \det C_{(n)}=Wr\bigl(\Phi_I(x)\bigr)|_{x=0} =\Delta_I.
$$
Now \eqref{3.34} implies
$$
\tau_I(t)=\frac1{\Delta_I} R\left(\frac1g\right) r(\chi_I) \tau_V(t).
$$
Noting that
$$
\left( C_{(n)}^{-1}A\right)^I= C^{-1}_{(n)}A^I
$$
and
$$
\det C_{(n)} =Wr\bigl(f_k(x)\bigr)|_{x=0}= \sum\det A^I\Delta_I
$$
completes the proof of \thref{3.8}.
\qed

\smallskip\noindent
Suppose that linear functionals $\chi_i$ can be defined on all $z^k$, e.g.\
$\chi_i$ are of the form:
\beq
\chi_i=\sum_{j\ge0}\alpha_{ij}\partial_z^j |_{z=\lambda_i}
\label{3.36}
\eeq
for $\lambda_i\ne0$ (cf.\ \cite{W}). We recall that (see e.g.\ \cite{I})
$$\ker\prod_{i=1}^r(L_V-\lambda_i^N)^{d_i} =\span\Big\{
\partial_z^{k_i}\Psi_V(x,z)\Big|_{z=\varepsilon^j\lambda_i}\Big\}_
{0\leq j \leq N-1, \,
1\leq i \leq r, \,
0\leq k_i \leq d_i-1} \, ,$$
with $\varepsilon=e^{2 \pi i/N},$ when  $\Psi_V(x,z)$ is well defined for
$z=\varepsilon^j\lambda_i$ (cf.\ eq.\ \eqref{1.13}).
Then we can put
\beq
f_k(t)=\langle c_k,\Psi_V(t,z)\rangle,\quad 0\le k\le n-1
\label{3.37}
\eeq
-- $f_k(t)$ can be thought as obtained from $f_k(x)$ by applying the KP flows.
In this case $\tau_W$ is given by the following theorem.
\bth{3.12}
If\/ $f_k(t)$ are as above and $g(z)=z^n$,
then
\beq
\tau_W(t)=\frac{Wr\bigl(f_k(t)\bigr)}{Wr\bigl(f_k(0)\bigr)}\tau_V(t).
\label{3.38}
\eeq
\eth
\noindent
{\it Proof\/} uses the {\it differential Fay identity\/} \cite{AvM2} (see
also \cite{vM}):
$$
Wr\bigl(\Psi_V(t,z_0),\ldots,\Psi_V(t,z_n)\bigr)\tau_V(t)
$$
$$
{}=\prod_{0\le j<i\le n} (z_i-z_j)\cdot\exp\Bigl(
\sum_{k=0}^{\infty}\sum_{i=0}^n t_kz^k_i\Bigr)
\tau_V\Bigl(t-\sum_{i=0}^n[z_i^{-1}]\Bigr),
$$
where $[z^{-1}]=(z^{-1}, z^{-2}/2, z^{-3}/3,\ldots)$.
After introducing the vertex operator
$$
X(t,z)=\exp\left(\sum_{k=1}^\infty t_k z^k\right)
\exp\left(-\sum_{k=1}^\infty \frac1{kz^k}\frac{\partial}{\partial t_k}\right)
$$
the RHS can be written in the form
$$
z_0^0z_1^1\cdots z_n^nX(t,z_n) X(t,z_{n-1}) \cdots X(t,z_0)\tau_V(t)
$$
$$
{}=z^n_n X(t,z_n)\Bigl(Wr\bigl(
\Psi_V(t,z_0),\ldots,\Psi_V(t,z_{n-1})\bigr)\tau_V(t)\Bigr).
$$
We apply the condition $c_0$ to the variable $z_0$, $c_1$ to $z_1$, $\ldots$,
$c_{n-1}$ to $z_{n-1}$, and set $z_n=z$. We obtain
$$
Wr\bigl(f_0(t),\ldots, f_{n-1}(t),\Psi_V(t,z)\bigr)\tau_V(t)
$$
$$
{}=z^nX(t,z) \Bigl( Wr\bigl(f_0(t),\ldots,f_{n-1}(t)\bigr)\tau_V(t)\Bigr).
$$
But \eqref{3.25} with $g(z)=z^n$ imply
$$
\Psi_W(t,z)=\frac{Wr\bigl( f_0(t),\ldots, f_{n-1}(t),\Psi_V(t,z)\bigr)}
{z^n Wr\bigl(f_0(t),\ldots, f_{n-1}(t)\bigr)}
$$
(KP flows applied to \eqref{3.25}).
Because the tau-function is determined from \eqref{1.7} up to a multiplication
by a constant, \eqref{3.38} follows (when $\tau_V(0)=\tau_W(0)=1$).
\qed
\bex{9.7}
Let us consider the simplest plane in the
Sato's Grassmannian $V= H_{+} = \span \{ z^i | i=0,1, \ldots \}$. Then
$$\psi_V(t,z)=\exp\sum t_k z^k,\quad L_V=\partial_x,\quad \tau_V(t)=1.$$
Every linear functional on $H_{+}$ is a linear combination of conditions of the
type
$$ e(k,\lambda)=\partial_z^k|_{z=\lambda} $$
and $h\left(L_{(0)}\right) =h(\partial_x)$ is an operator with constant
coefficients. The set of rational solutions of the KP hierarchy
\cite{Kr, SW} coincides with the set of B\"acklund--Darboux transformations
of $H_{+}$.
The formula \eqref{3.38} with $\tau_V=1$ is called a ``superposition law for
wobbly solitons'' (cf.\ \cite{W}, eq.\ (5.7)). The so called polynomial
solutions of KP \cite{DJKM, SW} correspond to the case $h(z)=z^d,$ i.e.\
all conditions are supported at 0.
\qed
\eex
Without any constraints on conditions $C$, there is a weaker version of
Theorem~\ref{t3.12}.
\bpr{3.13}
 For $g(z)=z^n$
$$
\tau_W(x)=\frac{Wr\bigl(f_k(x)\bigr)}{Wr\bigl(f_k(0)\bigr)}\tau_V(x),
$$
where $\tau_W(x)=\tau_W(x,0,0,\ldots)$.
\epr
\noindent
{\it Proof}. Formula \eqref{1.7} implies
\beq
\Psi_W(x,z)=e^{xz}\left(1-\partial_x\log \tau_W(x)z^{-1}+\cdots\right).
\label{3.42}
\eeq
On the other hand
$$
z^{-n}(\partial_x^n +p_1(x)\partial_x^{n-1}+\cdots+ p_0(x))\Psi_V(x,z)
$$
$$
{}=e^{xz}(1-\partial_x\log\tau_V(x) z^{-1} + p_1(x) z^{-1} +\cdots).
$$
Comparing the coefficients at $z^{-1}$ and noting that
$$
p_1(x)=-\partial_x\log Wr\bigl(f_k(x)\bigr)
$$
completes the proof.
\qed
%%%%%%%%%%%%%%%%% References %%%%%%%%%%%%%%%%%%%%%%%%%%%%%%%%%%%%%%%%%%%%%%%
\begin{small}
    
\end{small}
%%%%%%%%%%%%%%%%%%%%%%%%%%%%%%%%%%%%%%%%%%%%%%%%%%%%%%%%%%%%%%%%%%%%%%%%%%%%%%%
%%%%%%%%%%%%%%%%%%%%%%%%%%%%%%%%%%%%%%%%%%%%%%%%%%%%%%%%%%%%%%%%%%%%%%%%%%
\end{document}